\documentclass[12pt]{article}
\usepackage{amsxtra,amssymb,amsthm,amsmath,latexsym}

\theoremstyle{plain}

\newtheorem{theorem}{Theorem}[section]

\newtheorem{remark}{Remark}[section]
\numberwithin{equation}{section}

\def\R{{\mathbb R}}

\def\oH{{\overset{\circ}{H}}}
\def\oH1{{\overset{\circ}{H}\kern-.02in{}^1}}

\def\Im{{\hbox{\,Im\,}}}

\def\bee{\begin{equation*}}
\def\eee{\end{equation*}}
\def\be{\begin{equation}}
\def\ee{\end{equation}}

\begin{document}

\title{ A short proof of the existence of the solution to elliptic boundary problem}

\author{ A. G. Ramm
\\Mathematics Department, Kansas State University,\\
Manhattan, KS 66506-2602, USA\\
email: ramm@math.ksu.edu}

\date{}

\maketitle\thispagestyle{empty}

\begin{abstract}
There are several methods for proving the existence of the solution to the elliptic boundary problem
$Lu=f \text{\,\, in\,\,} D,\quad u|_S=0,\quad    (*)$. Here $L$ is an elliptic operator of second order,
$f$ is a given function, and uniqueness of the solution to problem (*) is assumed. The known methods
for proving the existence of the solution to (*) include variational methods, integral equation methods,
method of upper and lower solutions. In this paper a method based on functional analysis is proposed.
This method is conceptually simple and technically is easy. It requires some known a priori estimates and a continuation in a parameter
method.
\end{abstract}

\footnote{2010 Math subject classification: 35J05; 35J25}
\footnote{Key words:  elliptic boundary problems; continuation in a parameter.
}


\section{INTRODUCTION.}\label{Sec1}
Consider the boundary problem
\begin{equation}\label{eq:1}Lu=f   \text{ \,\, in \,\,} D,\end{equation}
\begin{equation}\label{eq:2}u=0   \text{\,\, on \,\,} S,\end{equation}
where $D\subset \R^3$ is a bounded domain with a $C^2-$smooth boundary $S$, $L$
is an elliptic operator,
\begin{equation}\label{eq:3}
Lu=-\partial_i(a_{ij}(x)\partial_j u) +q(x)u.
\end{equation}
Here and below $\partial_i=\frac {\partial}{\partial x_i}$, over the repeated indices summation is understood,
$1\le i,j\le 3,$ $a_{ij}(x)=a_{ji}(x),$ $\Im a_{ij}(x)=0,$
\begin{equation}\label{eq:4}
c_0|\xi|^2\le a_{ij}(x)\xi_i \overline{\xi_j}\le c_1|\xi|^2, \quad \forall x\in D,
\end{equation}
where $c_0, c_1 >0$ are constants independent of $x$ and $|\xi|^2=\sum_{j=1}^3 |\xi_j|^2$.
We assume that $q(x)$ is a real-valued bounded function and $|\nabla a_{ij}(x)|\le c$. One may easily
consider by our method the case of complex-valued $q$, see Remark 2.2 in Section 2.
By $c>0$ various estimation constants are denoted. In this paper the Hilbert space $H:=H^0:=L^2(D)$,
the Sobolev space $H^1_0$, the closure of $C^\infty_0(D)$ in the norm of the Sobolev space $H^1=H^1(D)$,
and the Sobolev space $H^2_0:=H^2(D)\cap H^1_0$ are used.

We assume for simplicity that problem (\ref{eq:1})-(\ref{eq:2}) has no more than one solution. This, for example,
is the case if
\begin{equation}\label{eq1.5}
(Lu,u)\ge c_2(u,u),   \quad \forall u\in D(L),
\end{equation}
where $c_2>0$ is a constant, and $D(L)=H^2_0$. The norm
in the Sobolev space $H^\ell$ is denoted  by the symbol $||\cdot||_\ell$. For example,
 \begin{equation}\label{eq1.6}
||\cdot||_2=\Big(\int_D(|u|^2+|\partial u|^2+|\partial^2u|^2)dx\Big)^{1/2}.
\end{equation}
By $|\partial u|^2$ the sum of the squares of the derivatives of the first order
is denoted, and  $|\partial^2 u|^2$ is understood similarly.

There is a large literature on elliptic boundary problems (see \cite{[A]}, \cite{[GT]}, \cite{[H]},  \cite{[LU]}, \cite{[L]}, \cite{[Mi]}) to name just a few books. Several methods were suggested to study problem \eqref{eq:1} - \eqref{eq:2}:  Hilbert space method, based on the Riesz theorem about bounded linear functionals (\cite{[LU]}, \cite{[L]}), integral equations of the potential theory (\cite{[Mi]}), method of lower and upper solutions (\cite{[GT]}).

The goal of this paper is to suggest a method for a proof of the existence of the solution to problem \eqref{eq:1} - \eqref{eq:2}, based on functional analysis. This method is simple, short, and does not require too much of a background knowledge from the reader.

The background material, that is used in our proof, includes the notions of closed  linear unbounded operators and symmetric operators (see \cite{[K]}) and second basic elliptic inequality (see \cite{[A]}, \cite{[GT]}, \cite{[LU]}, \cite{[L]}):
\begin{equation}\label{eq1.7}
||Lu||_0 \geq c_3||u||_{2}, \quad \forall u \in D(L),
\end{equation}
and the definition and basic properties of the mollification operator, see, for example, \cite{[A]}.

Let us outline the ideas of our proof. Let $R(L)$ denote the range of $L$.

We prove that  $R(L)$ is a closed subspace of $H^0$ and $R(L)^\perp = \{ 0 \}$.

This implies that $R(L) = H$, that is, problem \eqref{eq:1} - \eqref{eq:2} has a solution. Uniqueness of the solution follows trivially from the assumption \eqref{eq1.5}.

Let us summarize our result. This result is known (see, for example, \cite{[GT]}, \cite{[LU]}, \cite{[L]}), but we give a
short and essentially self-contained proof of it.

\begin{theorem}\label{Thm1.1}
Assume that $S$ is $C^2$- smooth, inequalities \eqref{eq:4}, \eqref{eq1.5} hold, and $q$ is a real-valued bounded function. Then problem \eqref{eq:1} - \eqref{eq:2} has a solution in $H_0^2$ for any $f \in H^0$, and this solution is unique. The operator $L$ is an isomorphism of $H^2_0(D)$ onto $H^0=L^2(D)$.
\end{theorem}

\begin{remark}\label{Rmk1.1}
We are not trying to formulate the result in its maximal generality. For example, one may consider by the same method
elliptic operators which are non-self-adjoint. In Section 2, Remark 2.2 addresses this question.
\end{remark}

In Section \ref{Sec2} proofs are given.

\section{Proofs}\label{Sec2}
It follows from \eqref{eq1.5} that
\begin{equation}\label{eq2.1}
||Lu|| \geq c_2||u|| \quad \forall u \in D(L), \quad ||u|| := ||u||_0.
\end{equation}
Therefore, if $Lu = 0$ then $u = 0$. This proves the uniqueness of the solution.

To prove the existence of the solution it is sufficient to prove that the range of $L$ is closed and its orthogonal complement
in $H^0$ is just the zero element. Indeed, one has
\begin{equation}\label{eq2.2}
H = \overline{R(L)} \bigoplus R(L)^\perp,
\end{equation}
where $R(L)^\perp$ denotes the orthogonal complement in $H = H^0$ and the over-line denotes the closure. Therefore, if
\begin{equation}\label{eq2.3}
R(L) = \overline{R(L)},
\end{equation}
and
\begin{equation}\label{eq2.3'}
R(L)^\perp = \{0\},
\end{equation}
then
\begin{equation}\label{eq2.4}
R(L) = H,
\end{equation}
and Theorem \ref{Thm1.1} is proved.

The closedness of $R(L)$ follows from inequality \eqref{eq1.7}. Indeed, if $Lu_n \xrightarrow[H^0]{} f$ then, by \eqref{eq1.7} and \eqref{eq2.1}, $u_n \xrightarrow[H_0^2]{} u$, so $u \in H^2_0(D):= D(L)$ and $Lu = f$. A more detailed argument goes as follows.
Let $v \in D(L)$ be arbitrary. Then
\begin{equation}\label{eq2.5}
(f, v) \xleftarrow[n \to \infty]{} (Lu_n, v) = (u_n, Lv) \xrightarrow[n \to \infty]{} (u, Lv), \quad \forall v \in D(L).
\end{equation}
Inequality \eqref{eq1.7} implies that $u \in H_0^2 = D(L)$. Therefore, formula \eqref{eq2.5} implies $Lu = f$. This argument proves that $R(L)$ is a closed subspace of $H^0$ and the operator $L$ is closed on $D(L)$.

Let us now prove that $R(L)^\perp = \{0\}$. Assume the contrary. Then there is an element $h \in H^0$ such that
\begin{equation}\label{eq2.6}
(Lu, h) = 0, \quad \forall u \in D(L) = H_0^2.
\end{equation}
 Let us derive from \eqref{eq2.6} that $h = 0$. To do this, first assume that $L =L_0:= -\Delta$, where $\Delta$ is the Dirichlet  Laplacian,
 and prove that $L_0$ is an isomorphism of $H^2_0(D)$ onto $H^0=L^2(D)$. This will prove Theorem 1.1 for $L=L_0$. Then we use continuation
 in a parameter method and prove that the same is true for $L$, which will prove Theorem 1.1.

 Take an arbitrary point $x \in D$, choose $\epsilon > 0$ so that the distance $d(x, S)$ from $x$ to $S$ is larger than $\epsilon$, and set $u = w_\epsilon(|x-y|)$, where $w_\epsilon(|x|)$ is a mollification kernel (see, for example, \cite{[A]}, p.5 ). This implies that $w_\epsilon(|x|) \in C_0^\infty(D) \subset D(L)$, and
\begin{equation}\label{eq2.7}
\lim_{\epsilon \downarrow 0}||\int_D w_\epsilon(|x - y|)h(y)dy - h(x)|| = \lim_{\epsilon \downarrow 0}||w_\epsilon * h - h|| = 0,
\end{equation}
where $w_\epsilon * h$ denotes the convolution. Then equation \eqref{eq2.6} yields
\begin{equation}\label{eq2.8}
-\int_D \Delta_y w_\epsilon(|x - y|)h(y)dy= - \Delta_x w_\epsilon * h = 0, \quad x \in D.
\end{equation}
Multiply \eqref{eq2.8} by $\eta_\epsilon:= w_\epsilon * h$, integrate over $D$, and then integrate by parts, taking into account that $\eta_\epsilon = 0$ on $S$ if dist$(x, S) > \epsilon$. The result is
\begin{equation}\label{eq2.9}
\int_{D}|\nabla \eta_\epsilon(x)|^2 dx = 0.
\end{equation}
From \eqref{eq2.9} it follows that $\nabla \eta_\epsilon=0$ in $D$, so $\eta_\epsilon=const$ in $D$.
Since this constant vanishes at the boundary $S$, it is equal to zero. Thus
\begin{equation}\label{eq2.10}
\eta_\epsilon(x) = w_\epsilon * h = 0 \quad \text{in} \quad D.
\end{equation}
Let $\epsilon \downarrow 0$ in \eqref{eq2.10} and get $h = 0$ in $D$. Thus, $R(\Delta)^\perp = \{0\}$, so $R(-\Delta) = H^0 = L^2(D)$.

Let us now prove that $R(L) = H^0$ for the operator \eqref{eq:3}. This is proved by a continuation in a parameter.
 Define $L_s =L_0+s(L-L_0)$, $0 \leq s \leq 1$, $L_0 = -\Delta$, $L_1 = L$. We prove that $R(L_s) = H^0$ for all $s\in [0,1] $
 and the map $L_s: H^2_0\to H^0$ is an isomorphism. For $s = 0$ this is proved above.

Consider equation \eqref{eq:1} with $L = L_s$ and apply the operator $L_0^{-1}$ to this equation. The result is
\begin{equation}\label{eq2.13}
u + sL_0^{-1}(L-L_0)u = L_0^{-1}f.
\end{equation}
This equation is in the space $H_0^2$. The norm of the operator $ sL_0^{-1}(L-L_0)$ in $H_0^2$ is less than one if $s$ is sufficiently small. Indeed, inequality similar to \eqref{eq1.7} holds for $L_s$ for any $s\in [0,1]$ with the same constant $c_3$, because this constant depends only on the bounds on the coefficients of $L_s$ and these bounds are independent of $s \in [0,1]$. Thus,
\begin{equation}\label{eq2.14}
||L_s u||_0 \geq c_3||u||_{2}, \quad \forall u \in H_0^2, \quad 0 \leq s \leq 1.
\end{equation}
Therefore,
\begin{equation}\label{eq2.14'}
||L_0^{-1}(L-L_0)u||_{2} \leq \frac{1}{c_3}||(L-L_0)u||_0 \leq c_3'||u||_{2}, \qquad \forall u \in H_0^2,
\end{equation}
because $||(L-L_0)u||_0\le c||u||_2$, where $c>0$ is a constant not depending on $s$, $c$ depends only on the bounds
on the coefficients of $L$.
Consequently, if $sc_3' < 1$, that is, if $s < (c_3')^{-1}$, then equation \eqref{eq2.13} is
 uniquely solvable in $H_0^2$ for any $f \in H^0$, and $R(L_s) = H^0$.

 Let $s_0:= \frac{1}{2} (c_3')^{-1}$. Define $L_{s_0+s'} =L_{s_0} + s'(L-L_{s_0})$, $s'\in [0,1]$. One has  $L_{s_0+s'} =L_{s_0}$ as $s'=0$
 and  $L_{s_0+s'} =L$ as $s'=1$.  Applying the same argument and
 using the fact that $||L^{-1}_{s_0}||_{H^0\to H^2_0(D)}$ does not depend on $s_0$,   one gets
\begin{equation}\label{eq2.15}
||L^{-1}_{s_0}(L-L_{s_0})u||_{2} \leq c_3'||u||_{2}.
\end{equation}
Therefore, for $s' < (c_3')^{-1}$, one has
\begin{equation}\label{eq2.16}
||s'L_{s_0}^{-1}(L-L_{s_0})|| < 1.
\end{equation}
Let $s'=s_1:=\frac 1 2  (c_3')^{-1}$. Then $R(L_s)=H^0$ for $s<s_0+s_1$ and $L_s: H^2_0\to H^0$  is an isomorphism.
Consequently, repeating the above argument finitely many times one reaches the operator $L$ and gets both conclusions:  $R(L) = H^0$ and $L$ is an isomorphism of $H_0^2$ onto $H^0$.

Theorem \ref{Thm1.1} is proved. \hfill$\Box$

\begin{remark}\label{Rmk2.1}
The method of continuation in a parameter goes back to \cite{[S]}, see also \cite{[L]}.
\end{remark}

\begin{remark}\label{Rmk2.2}
Consider the operator $L_1 = L + L'$, where $L'$ is an arbitrary first order differential operator and $L$ is the same as in Section \ref{Sec2}. The operator $L_1$ is not necessarily symmetric. Problem \eqref{eq:1} - \eqref{eq:2} is equivalent to the operator equation
\begin{equation}\label{eq2.11}
u + Au = L^{-1}f \quad \text{in} \quad H^0,
\end{equation}
where
\begin{equation}\label{eq2.12}
A = L^{-1}L'
\end{equation}
is a compact operator in $H^0$. This follows from the Sobolev embedding theorem (\cite{[A]}, \cite{[GT]}).

Therefore, the Fredholm alternative holds for equation \eqref{eq2.11}. So, if the homogeneous version of the equation \eqref{eq2.11} has only the trivial solution (zero solution) then equation \eqref{eq2.11} is solvable for any $f$, and its solution $u \in H_0^2$.
\end{remark}

\begin{remark}\label{Rmk2.3}
If $L$ is symmetric on $D(L)=H^2_0(D)$, then Theorem 1.1 shows that $L$ is self-adjoint on $D(L)$. Indeed, the definition
of the adjoint operator $L^*$ says, that $(Lu,w)=(u,w^*)$ for all $u\in D(L)$. By Theorem 1.1 there exists
$z\in D(L)$ such that $Lz=w^*$. Thus, $(Lu,w)=(u, Lz)=(Lu,z)$. Since the range $R(L)=H^0$, it follows that $w=z$. So,
$w\in D(L)$, $D(L^*)=D(L)$ and $L=L^*$, as claimed.
\end{remark}

\end{document}